\documentclass[10pt,twoside]{amsart}

\usepackage{amscd}
\usepackage{amssymb} 
\usepackage{verbatim}
\usepackage{amsmath}
\usepackage{amsxtra}
\usepackage[noadjust]{cite}         

\newtheorem{theorem}{Theorem}[section]
\newtheorem{lemma}[theorem]{Lemma}

\theoremstyle{definition}

\numberwithin{equation}{section}

\usepackage{amssymb}
\usepackage{amsfonts}
\usepackage{color}
\usepackage{graphics}
\usepackage{epsfig,psfrag,graphicx}
\usepackage{amssymb}
\usepackage{eepic,epic}
\usepackage{color}

\usepackage{epsfig} 
\usepackage{amsmath} 
\usepackage{amssymb} 
\usepackage{amsbsy} 

\newcommand{\nc}{\newcommand}
\nc{\R}{{\mathbb R}} 
\nc{\N}{{\mathbb N}}
\nc{\weak}{\rightharpoonup}
\nc{\weakstar}{\stackrel{\ast}{\rightharpoonup}} 
\renewcommand{\div}{{{\mathrm{div}}_x}\,}
\newcommand{\vrho}{\varrho}
\nc{\modular}[1]{{\stackrel{ #1}{\longrightarrow\,}}}

\def\bbbone{{\mathchoice {\rm 1\mskip-4mu l}
{\rm 1\mskip-4mu l} {\rm 1\mskip-4.5mu l} {\rm 1\mskip-5mu l}}}

\def\tens#1{\pmb{\mathsf{#1}}}
\def\vec#1{\boldsymbol{#1}}

\newcommand{\vr}		{\vrho}
\newcommand{\vre}		{\vr_\varepsilon}
\newcommand{\vrez}	{\vr_{0,\varepsilon}}
\newcommand{\vrs} 		{\overline{\vr}}
\newcommand{\vret}		{\tilde{\vr}_\ep}
\newcommand{\ue}		{\vec{u}_\varepsilon}
\newcommand{\uez}		{\vec{u}_{0,\varepsilon}}
\newcommand{\ep}		{\varepsilon}
\newcommand{\temp}	{\vartheta}
\newcommand{\tem}		{\vartheta_\varepsilon}
\newcommand{\temz}	{\vartheta_{0,\varepsilon}}
\newcommand{\tems}	{\overline{\vartheta}}
\newcommand{\tS}		{\tens{S}}
\newcommand{\q}		{\vec{q}}
\newcommand{\n}		{\vec{n}}
\newcommand{\vU}		{\vec{U}}
\newcommand{\ess} 	{{\rm{ess}}}
\newcommand{\res}		{{\rm{res}}}
\newcommand{\dx}		{\,{\rm d}x}
\newcommand{\dt}		{\, {\rm d}t}
\newcommand{\dxdt}	{\, {\rm d}x{\rm d}t}
\newcommand{\len}		{\Delta_{\ep, {\rm N}}}

\newcommand{\U}		{\vec{U}}

\def\bbbone{{\mathchoice {\rm 1\mskip-4mu l}
{\rm 1\mskip-4mu l} {\rm 1\mskip-4.5mu l} {\rm 1\mskip-5mu l}}}

\date{}

\begin{document}




\title[The low Mach number limit]{Stability with respect to domain of the low Mach number limit of compressible heat-conducting viscous fluid}

\author[A. Wr\'oblewska-Kami\'nska]{Aneta Wr\'oblewska-Kami\'nska}
\address{Institute of Mathematics of the Polish Academy of Sciences, ul. \'Sniadeckich 8, 00-656 Warszawa, Poland} 
  \email{email: awrob@impan.pl}

\thanks{This research was supported by the grant Sonata Bis UMO-2020/38/E/ST1/00469, National Science Centre, Poland}


%
%
%
%
%
%
%
\begin{abstract}
We investigate the asymptotic limit of solutions to the Navier--Stokes--Fourier system with the Mach number proportional to a small parameter $\ep \to 0$, the Froude number proportional to $\sqrt{\ep}$ and when the fluid occupies large domain with spatial obstacle of rough surface varying when $\ep \to 0$.  The limit velocity field is solenoidal and satisfies the incompressible Oberbeck--Boussinesq approximation. Our studies are based on weak solutions approach and
in order to pass to the limit in a convective term we apply the spectral analysis of the associated wave propagator (Neumann Laplacian) governing the motion of acoustic waves. 
\end{abstract}





\maketitle

\section{Introduction and formulation of the problem}

The Oberbeck--Boussinesq  approximation is a mathematical model of a stratified flow, where the fluid is assumed to be incompressible and yet convecting a diffusive quantity creating positive and negative buoyancy force. Then the system of equations reads:
	\begin{equation}\label{limit11}\tag{OB$^1$}
	\div \vU =0,
	\end{equation}
	\begin{equation}\label{limit12}\tag{OB$^2$}
	\vrs\left(\partial_t \vU + \div (\vU \otimes \vU)  \right) + \nabla_x P 
	= \mu \Delta \vU + r \nabla_x F,
	\end{equation}
	\begin{equation}\label{limit13}\tag{OB$^3$}
	\vrs c_p\left( \partial_t \Theta + \div (\vU \Theta) \right) - \kappa(\tems)\Delta \Theta - \vrs\tems\alpha \div(F\vU)= 0, 
	\end{equation}
	\begin{equation}\label{limit14}\tag{OB$^4$}
	r + \vrs \alpha \Theta = 0,
	\end{equation}
where $\vU$ denotes the velocity of the fluid, $\Theta$ stands for the deviation of the temperature, 
$P$ is the pressure, constants $\mu,$ $\kappa,$ $\vrs,$ $c_p$, $\alpha$ are positive (will be defined later).
Here $F$ stands for potential of a driving force (e.g. gravitational potential) acting on the fluid. Let us note that the density is constant in the Oberbeck--Boussinesq approximation except in the buoyancy  force, where it is interrelated in the temperature deviation through Boussinesq relation \eqref{limit14}, (see Zeytounian \cite{Ze2003}). Let us notice that in the OB approximation 
$\Theta$ is a deviation of temperature from the equilibrium rather then the temperature itself
and the temperature deferences are not caused by the flow, but exists independent of the flow.
Our aim is to derive the above system on an exterior domain $\R^3 \setminus O$ with no-slip boundary condition on the bounded obstacle $O$. 
Therefore we study stability of the  rescaled compressible  Navier--Stokes--Fourier system when a Mach number is proportional to 
a small parameter, i.e. $Ma=\ep$ and $\ep \to 0$, and a Froude number $Fr = \sqrt{\ep}$.  About other characteristic numbers like Strouhal, Reynolds, P\'eclet number we assume  they are equal one.

We are motivated by a similar asymptotic analysis of barotropic compressible fluid flow, described by the Navier--Stokes system with a low Mach number on varying domains provided in \cite{FKKS}. Our aim is to extend this  result to the case of heat-conducting fluids by methods developed in \cite{F_LDB,FS}. The asymptotic analysis of complete fluid system on varying domains (but in different way then here) and with a small Mach number is considered in \cite{AWK}, where the author justify OB system on whole $\R^3$ space with concentric gravitation force.

Following \cite{BFNW,FKKS} we introduce a class of admissible domains with rough (oscillating) boundaries of some obstacle. It was observed that such a choice may give rise to the no-slip boundary condition for the asymptotic limit of velocity field.
In particular we assume that the given family of domains $\{ \Omega_\ep \}_\ep$ satisfies the following hypothesis:
\begin{description}
\item[D1.)] $\Omega_\ep \subset \R^3$ is bounded domain with $C^2$ boundary for each $\ep \in (0,1)$ and $\partial \Omega_\ep = \partial O_\ep  \cup \mathcal{S}_{\ep}$; 
\item[D2.)] for simplicity we assume that the outer part of boundary $\mathcal{S}_{\ep}$ consists of a sphere centred in the origin and of a radius $\frac{1}{\ep^\delta}$ with  $\delta > 0$ (i.e. the domain is sufficiently ''large'');
 \item[D3.)] the boundary of the obstacle $\partial O_\ep$ is such that for all $\ep \in (0,1)$ $O_\ep  \subset B_r(0) \subset B_{1/\ep^\delta} (0)$ with some fixed $r >0$;
  \item[D4.)] $\R^3 \setminus O_\ep$ satisfies the uniform $\alpha$-cone condition with $\alpha >0$ independent of $\ep$. Namely for any $x_0 \in \partial O_\ep$ there exists a unit vector $\xi_{x_0} \in \R^3$  s.t. 
 $C(x,\gamma,\alpha, \xi_{x_0}) \subset (\R^3\setminus O_\ep)$ whenever $x \in \R^3 \setminus O_\ep,$ $|x-x_0| < \alpha$, where 
 $C(x,\gamma,\alpha, \xi) = \{ y \in \R^3\,| \, 0 < | y-x| \leq \alpha, (y - x ) \cdot \xi > \cos (\gamma) |y-x| \}$ with vertex at $x$, aperture $2\gamma < \pi$, height $\alpha$, and orientation given by a unit vector $\xi$;
 \item[D5.)] for each $x_0 \in \partial O_\ep,$ there are two open balls $B_r(x_i) \subset \Omega_\ep$, $B_r(x_j) \subset O_\ep$ of radius  $r > c_b\ep^\beta$ (the radius $r$ may change but sufficiently ''slow'') such that 
     	$ \overline{B_r(x_i)} \cap \overline{B_r (x_j)} = x_0$  with $c_b >0,$ $\beta > 0$  independent of $\ep$;
\item[D6.)] after translation and rotation of the coordinate system, a part $\Gamma \subset \partial O$ can be described by a graph of function $\gamma \in W^{1,\infty}(U),$ $U \subset \R^3$ and $\Gamma = \{ x\in \R^3\ : \  (x_1,x_2)\in U, \ x_3 = \gamma(x_1,x_2)\}$ while $\Gamma_\ep = \partial O_\ep \cap U \times \R$ are represented 
by $\Gamma_\ep= \{ x\in \R^3\ : \  (x_1,x_2)\in U, \ x_3 = \gamma_\ep(x_1,x_2)\}$, where $\{ \gamma_\ep \}_{\ep}$ is a bounded sequence in $W^{1,\infty}(U)$, $\gamma_\ep \to \gamma$ in $C(\overline{U})$. Moreover $\Gamma_\ep$ are oscillating for $\ep \to 0$. Namely, when we introduce a Young measure $\nu[y]$, $y\in U$, associated to the sequence 
$\{ \nabla_y \gamma_\ep\}_{\ep}$, we suppose that ${\rm supp\,[\nu[y]]}$ contains two independent vectors in $\R^2$ for a.a. $y \in U.$  
\end{description}
In certain sense $\Omega_\ep \to \R^3\setminus O.$
 We give here a mathematical justification of the Oberbeck--Boussinesq approximation of a incompressible flow on exterior domain $\Omega  = \R^3 \setminus O$ with no-slip boundary condition on the obstacle by asymptotic analysis of weak solutions to the compressible  Navier--Stokes--Fourier system in a low  Mach number regime: $Ma=\ep \to 0$, on a family of domain $\Omega_\ep$ varying with $\ep >0$.

\section{Primitive system}

In the beginning of this section let us introduce some standard notation.
We denote by $\langle \cdot ,  \cdot \rangle $ duality pairing.
By $L^p(B)$ we mean the space of Lebesgue measurable functions $g$, where $|g|^p$ is integrable over set $B$.
The Sobolev space of functions which derivatives are integrable up to order $k$ in $L^p$ we denote by $W^{k,p}$.
By $\mathcal{D}^{k,p}(B)$ we set homogenous Sobolev spaces i.e.
$\mathcal{D}^{k,p}(B) = \{ g \in L^1_{\rm loc}(B)\, : \, D^\alpha g \in L^p(B),\ |\alpha| = k \}$, where 
$k \geq 0$ and $p \geq 1$. In the whole paper $c$ will denote generic constant which may change from line to line.
 \bigskip
 
We start our considerations with a "primitive system" - the rescaled Navier--Stokes--Fourier system with a small Mach and Froude number which consists of: the continuity equation (conservation of mass), the momentum equation, the entropy balance and the total energy balance respectively
	\begin{equation}\label{ceq}\tag{NSF$_{\ep}^1$}
	\partial_t \vre + \div (\vre\ue)=0,
	\end{equation}
	\begin{equation}\label{meq}\tag{NSF$_{\ep}^2$}
	\partial_t (\vre\ue)+ \div(\vre\ue\otimes\ue) +\frac{1}{\ep^2} \nabla_x p(\vre,\tem) 
	= \div \tS(\tem,\nabla_x\ue) + \frac{1}{\ep} \vre\nabla_x F_\ep ,
	\end{equation}
	\begin{equation}\label{eiq}\tag{NSF$_{\ep}^3$}
	 \partial_t (\vre s(\vre, \tem))  + \div (\vre s (\vre,\tem)\ue) + \div\left(\frac{\q(\tem,\nabla_x \tem )}{\tem} \right)
	= \sigma_\ep,
	\end{equation}
	\begin{equation}\label{eeq}\tag{NSF$_{\ep}^4$}
	\frac{d}{dt} \int_{\Omega_\ep} 
	\left( \frac{1}{2} \vre|\ue|^2 + \frac{1}{\ep^2} \vre e(\vre,\tem) - \frac{1}{\ep} \vre F_\ep \right) dx = 0.
	\end{equation}
Where the viscous stress tensor satisfies the Newton rheological law and heat flux is determined by the Fourier law:
	\begin{equation*}  
	\tS(\tem,\nabla_x \ue) = \mu(\tem)\left( \nabla_x\ue + \nabla_x^T \ue  - \frac{2}{3}\div \ue \tens{Id} \right) + \eta(\tem) \div\ue \tens{Id}
	\end{equation*}
	\begin{equation*}
	 \q(\tem,\nabla_x \tem)= - \kappa(\tem) \nabla_x \tem 
	\end{equation*}
with a positive heat coefficient $\kappa$ and for the entropy production rate holds: 
	\begin{equation}\label{ss}
	\quad
	 \sigma_\ep \geq \frac{1}{\tem} \left({\ep^2} \tS_\ep(\tem,\nabla_x\ue) : \nabla_x \ue
	 - \frac{\q_\ep(\tem,\nabla_x \tem )\cdot \nabla_x \tem}{\tem}  \right).
	\end{equation}
 The unknowns are the fluid mass density $\vre=\vre(t,x)$, the velocity field $\ue=\ue(t,x) : (0,T)\times \Omega_\ep \to \R^3$ 
and absolute temperature $\tem=\tem(t,x) : (0,T) \times \Omega_\ep \to \R$. The pressure $p$, the specific internal energy 
$e$ and the specific entropy $s$ are given scalar valued functions of $\vr$ and $\temp$ which are related through 
Gibbs' equation	
	$
	\temp D s = D e + p D ( {1}/{\vr}).	
	$
The system is supplemented with complete slip boundary conditions for velocity field and  the boundary of physical space is thermally isolated, i.e.
	\begin{equation}\label{bc1-2}
	\ue \cdot \n |_{\partial \Omega_\ep} = 0,
	\quad
	[ \tS (\tem, \nabla_x \ue) \n ] \times \n = 0, \quad\quad
	\q\cdot\n|_{\partial {\Omega_\ep}}=0 .
	\end{equation}		
{Small parameter $\ep$ in the system (\ref{ceq} - \ref{eeq}) results from  dimensionless form of a Navier--Stokes--Fourier system and corresponds to small Mach and Froude number (Ma$=\ep$, Fr$= \sqrt{\ep}$), see \cite{FN}, Klein at al. \cite{K_01}, Zeytounian \cite{Ze2004}. Smallness of Mach number physically means that characteristic speed of the flow is dominated by the speed of the sound in the medium under consideration.}  Assumption that ${\rm Fr} >> { \rm Ma}$ means that external sources of mechanical energy are small and  $\frac{\rm Ma}{ \rm Fr}\to 0$,  
what corresponds to low stratification).

 \subsection{Structural restrictions}
In order to be able to use the existence result of \cite{FN} and later to build uniform estimates,
we need to  impose structural restrictions on the thermodynamical functions $p$, $e$, $s$ as well as on  
the transport coefficients $\mu$, $\eta$, $\kappa$. Following \cite{FN} (where the reader can find more detailed description and physical motivations) we set
	\begin{equation}\label{pp1}
	\begin{split}
	& p(\vre,\tem)= \tem^{5/2} P\left(\frac{\vre}{\tem^{3/2}}\right) + \frac{a}{3} \tem^{4},\ a>0, 
	\mbox{ where }P\in C^1 [0,\infty)\cap C^2(0,\infty),\\
	& P(0)=0,\quad P'(Z)>0\mbox{ for all }Z\geq 0,
	\end{split}
	\end{equation}
	\begin{equation}\label{pp3}
	\begin{split}
	& 0<\frac{ \frac{5}{3} P(Z) - P'(Z)Z }{ Z }< c \mbox{ for all } Z > 0, \quad \lim\limits_{Z\to \infty} \frac{P(Z)}{Z^{5/3}} = P_\infty >0,  
	\\ & \mbox{ and } \partial_\vr  p (\vr,\temp)>0.
	\end{split}
	\end{equation}
Accordingly to Gibbs' relation, the specific internal energy and the entropy can be written in the following forms	
	\begin{equation}\label{int_en}
	e(\vr,\temp) = \frac{3}{2} \frac{\temp^{5/2}}{\vr} P\left( \frac{\vr}{\temp^{3/2}} \right)
	+ a\frac{\temp^{4}}{\vr}, \quad 	\partial_\temp e(\vr,\temp) >0, \quad \mbox{ is positive and bounded,}
	\end{equation}
	\begin{equation}\label{entro}
	s(\vr,\temp)= S\left(\frac{\vr}{\temp^{3/2}}\right) + \frac{4}{3} a \frac{\temp^{3}}{\vr},\quad\quad
	S'(Z) = -\frac{3}{2} \frac{ \frac{5}{3} P(Z) - Z P'(Z)}{Z^2} \mbox{ for all } Z>0.
	\end{equation}
The transport coefficients: $\mu$ - shear viscosity, $\eta$ - bulk viscosity and $\kappa$ - heat conductivity  are assumed to be continuously differentiable 
functions of the temperature  $\temp \in [0,\infty)$ satisfying the following  growth conditions for all $\temp\geq 0$ and some positive constants $\underline\mu$, $\overline\mu$, $\overline\eta$, $\underline\kappa$, $\overline\kappa$:
	\begin{equation}\label{mu}
	\begin{split}
	&0<\underline\mu(1+\temp) \leq \mu(\temp) \leq \underline\mu (1 + \temp),
	\quad\quad
	0 \leq \eta(\temp) \leq \overline\eta(1 + \temp), 
\\ &
	0 <  \underline\kappa (1 + \temp^3) \leq \kappa(\temp) \leq \overline\kappa(1+\temp^3).
	\end{split}
	\end{equation}

\subsection{Equilibrium state and ill-prepered initial data}

Let us assume that outer force $F$ is defined on whole space $\R^3$ and is independent of $\ep$. The so-called equilibrium state (static state) for each scaled NSF$_\ep$ system  consist of static density $\vret$ and constant temperature distribution  $\tems$
satisfying (for a convenience we consider a static density  $\vret$ defined on the whole space $\R^3$)
	\begin{equation*}\label{prF}
	\nabla_x p(\vret,\tems) = \ep \vret \nabla_x F_\ep  \mbox{ in }\R^3
	\quad
	\mbox{ where }\lim\limits_{|x|\to \infty} \vret(x) = \vrs.
	\end{equation*}
Hence we have 
	\begin{equation}\label{nff}
	\begin{split}
	& \vret - \vrs = \frac{\ep}{P'(\vrs)} F  + \ep^2 h_\ep F_\ep, 
	\mbox{ with } P'(\vr) = \frac{1}{\vr} \partial_\vr p(\vr,\tems),\  \| h_\ep \|_{L^{\infty}(\R^3)} < c 
	\\ & \mbox{ and }
	|\nabla_x \vret(x)| \leq \ep c | \nabla_x F_\ep (x) | \mbox{ for } x\in \R^3 
	\end{split}
	\end{equation}		
(notice that the above properties gives closeness of static density $\vret$ and constant state $\vrs$).
Since we work with weak solutions based on energy estimates and control of entropy production rate we need to assume that initial data are close to equilibrium state.
Namely initial density and initial temperature are of the following form
	\begin{equation}\label{in_vr}
	\vrez = \vret + \ep \vrez^{(1)}, \quad\quad
	\temz = \tems + \ep \temz^{(1)}
	\end{equation}
where
$\tems>0$  is positive constants characterising the static distribution 
of the absolute temperature and
	\begin{equation}\label{illdatabaud}
	\begin{split}
	&  \|  \vrez^{(1)} \|_{L^\infty \cap L^2(\Omega_\ep)} \leq c, \ \int  \vrez^{(1)} \dx = 0,\  \|  \temz^{(1)} \|_{L^\infty \cap L^2 (\Omega_\ep)} \leq c, \  \int  \temz^{(1)} \dx = 0,
	  \\ &
	  \  \| \vec{u}_{0,\ep}  \|_{L^\infty \cap L^2 (\Omega_\ep)} \leq c  \mbox{ for all }\ep \in (0,1].
	  \end{split}
	\end{equation}
The above uniform bounds will allow to control right hand side of total dissipation balance which is a source of uniform estimates needed to perform the limit system.
Nevertheless, such a choice allow to consider nontrivial dynamics  but on the other hand it causes  oscillations in acoustic equation. Those will  be eliminated by dispersive estimates.

\subsection{Main result}
We say that functions $\vU$, $\Theta$ and $r$ are a weak solution to the Oberbeck--Boussinesq approximation (OB) if holds:
$\vU\in L^\infty(0,T;L^2(\Omega;\R^3)) \cap L^2(0,T; W^{1,2}(\Omega;\R^3)),$  $\Theta \in L^\infty(0,T; L^2(\Omega)) \cap L^2(0,T; W^{1,2}(\Omega))$, 
$r \in L^\infty (0,T; L^{5/3}_{\rm loc} (\Omega))$ and 
	\begin{equation}\label{ob11}
	\begin{split}
	\div \vU & = 0 \mbox{ a.e. on }(0,T)\times \Omega,
	\\
	\int_0^T \int_{\Omega} (\vrs (\vU \cdot \partial_t \varphi + (\vU \otimes \vU) : \nabla_x \varphi )) \dx\dt
	& = \\
	-\int_{\Omega}   \vrs \vU_0 \cdot \varphi(0,\cdot) \dx 
	 + & \int_0^T \int_{\Omega}\left( \tS : \nabla_x \varphi - r \nabla_x F \cdot \varphi \right) \dx\dt
	\end{split}
	\end{equation}
for any {  $\varphi \in C^\infty_c( [0,T); {C^\infty_c}(\Omega; \R^3))$}, where $\div \varphi = 0$ and
	$\tS = \mu(\tems)(\nabla_x \vU + \nabla_x \vU^T).$
Moreover
	\begin{equation}\label{ob12}
	\begin{split}
	\vrs c_p (\vrs,\tems) \left[ \partial_t \Theta + \div (\Theta \vU )\right]
	&  - \div( \kappa(\tems) \nabla_x \Theta) 
	 - \vrs \tems \alpha(\vrs,\tems) \div(F\vU)
	 = 0 \\ & \mbox{ a.e. in } (0,T)\times \Omega
	\\
	\Theta(0,\cdot) 
	& = \Theta_0
	\\
	 r+ \vrs\alpha(\vrs,\tems)\Theta 
	& = 0  \mbox{ a.e. in } (0,T)\times \Omega. 
	\end{split}
	\end{equation}
By $c_p$ we mean specific heat at constant pressure 
and $c_p(\vrs, \tems) = \partial_\temp e(\vrs,\tems) +  \alpha(\vrs, \tems) \frac{\tems}{\vrs} \partial_\temp p (\vrs,\tems)$
by $\alpha > 0$ we mean the coefficient of thermal 
expansion of the fluid, $\alpha(\vrs, \tems) = \frac{1}{\vrs}   \frac{\partial_\temp p(\vrs,\tems)}{\partial_\vr p (\vrs,\tems)}$, 
both are evaluated at the reference density $\vrs$ and temperature $\tems$. Then the main result reads as follows:
	
\begin{theorem}\label{thmain}
Let $\Omega_\ep\subset \R^3$ be a family of domains defined by (D1)--(D5) with $\beta < \frac{1}{4}$ and $\delta > 1$. Assume that $p$, $e$, and $s$ satisfy (\ref{pp1}--\ref{entro}), the transport coefficients $\mu$, $\eta$ and $\kappa$ satisfy growth conditions  \eqref{mu} and driving force is determined by a scalar potential $F \in W^{1,\infty} (\R^3)$.
Let $\{\vre,\ue,\tem\}_{\ep>0}$ be a family of weak solutions to the scaled Navier--Stokes--Fourier system (\ref{ceq}--\ref{eeq}), on the sets $(0,T) \times \Omega_\ep$, supplemented with boundary conditions (\ref{bc1-2}) and initial data 
(\ref{in_vr}) with $\vret >0$, $\vrs >0$ and $\tems >0,$ and  satisfying (\ref{illdatabaud}) for all $\ep \in (0,1)$.
Moreover we assume that 
	\begin{equation*}\label{asiv}
	\begin{split}
	& \vrez^{(1)} \weak\vr_0^{(1)} \mbox{ weakly in }L^2(\R^3), 
	 \quad 
	 { \uez \weak \vU_0 \mbox{ weakly in }L^2(\R^3;\R^3), }
	\\ & 
	 \temz^{(1)} \weak \temp_0^{(1)} \mbox{ weakly in } L^2(\R^3).
	\end{split}
	\end{equation*}
	Then for suitable subsequence as $\ep \to 0$ we obtain that
	\begin{equation*}
	\begin{split}
	& \vre \to \vrs \mbox{ strongly in } L^\infty(0,T; L^{5/3}(K)),
	\quad \quad 
	\frac{\vre - \vrs }{\ep} \weak r \mbox{ weakly in }L^2(0,T; L^2 (K)),
	\\ &
	\frac{\tem - \tems}{\ep} \weak \Theta \mbox{ weakly in }L^2(0,T; W^{1,2}(\R^3)),
	\\ &
	\ue \weak \vU \mbox{ weakly in } L^2(0,T; W^{1,2}(\R^3; \R^3)),
	\ 
	\ue \to \vU \mbox{ strongly in } L^2((0,T) \times K; \R^3 )
	\end{split}
	\end{equation*}
for any compact set { $K \subset \Omega $}, where functions $\vU$, $\Theta$ is a weak solution of the Oberbeck--Boussinesq approximation (\ref{limit11}--\ref{limit14}) in $(0,T) \times \R^3$  in the sense specified in (\ref{ob11}--\ref{ob12})
	with $\vU(0,\cdot) = \vec{H}[\vU_0]$ and $\Theta_0 = \temp_0^{(1)}.$ Moreover if (D6) is satisfied, $\vU|_{\partial O} = 0.$
\end{theorem}
Here $\vec{H}[ \cdot]$ denotes the projection on the space of divergence free functions on $\Omega$  of Helmholtz decomposition.
 The rest of the paper is devoted to the proof of the Theorem~\ref{thmain} or rather to the sketch of the proof with references where reader can find all details.  

\section{Proof of the Theorem~\ref{thmain}}

Since for each $\ep \in (0,1)$  the set $\Omega_\ep$ is sufficiently regular and bounded, in order to provide the existence of the family of weak solutions  $\{\vre,\ue,\tem\}_{\ep>0}$ to  the primitive system - compressible Navier--Stokes--Fourier (\ref{ceq} - \ref{eeq}) stated on $\Omega_\ep$ we use the result of E. Feireisl and A. Novotn\'y \cite[Theorem 3.2]{FN}. Then the  following regularity of solutions  can be obtained: 
	$\vre \in C_{{\rm weak}}( 0,T; L^{5/3}(\Omega_\ep))$, 
	 $\vre \in L^q ((0,T)\times \Omega_\ep)$ for a certain $q>\frac{5}{3}$ and 
	$\ue \in L^2(0,T; W^{1,2}(\Omega_\ep;\R^3)).$ Moreover 
the absolute temperature $\tem$ is a measurable function $\tem(t,x) >0$ for a.a. $(t,x)\in (0,T) \times \Omega_\ep$ 	and 
	$\tem \in L^2(0,T; W^{1,2}(\Omega_\ep))\cap L^{\infty} (0,T; L^4 (\Omega_\ep)),$
	$\log \tem \in L^2(0,T; W^{1,2}(\Omega_\ep)).$
	
\subsection{Uniform bounds}\label{unib}
All uniform bounds stated below may be seen as a direct consequence of total dissipation balance and more detailed reasoning may be found in \cite{FN,FS,AWK}.

To begin with, according to these references, we introduce essential and residual part of  a measurable function $h$ as
	\begin{equation*}\label{ess-def}
	h = [h]_{\ess} + [h]_{\res},\  [ h]_{\ess} = \chi(\vre,\tem) h, \  [h]_{\res} = (1-\chi(\vre,\tem))h,
	\end{equation*}
where $\chi \in C^\infty_c ( (0,\infty)\times(0,\infty) )$, $ 0\leq \chi \leq 1$, $\chi=1$ on the set  ${\mathcal{O}}_{\ess}$ and $
	 {\mathcal{O}}_{\ess} =  [\vrs/2 , 2 \vrs ] \times [\tems/2, 2 \tems],$
	$ {\mathcal{O}}_{\res} = (0,\infty)^2\setminus {\mathcal{O}}_{\ess} .$
%
The total dissipation balance reads then
\begin{equation}\label{H2}
	\begin{split}
	& \int_{\Omega_\ep} \left( \frac{1}{2}\vre|\ue|^2 \right)(t) \dx  +\frac{1}{\ep^2}\left(H_{\tems}(\vre,\tem) - (\vre - \vret)\frac{\partial H_{\tems}(\vret,\tems)}{\partial \vr} 
	- H_{\tems}(\vret,\tems)\right)(t) \dx
	\\ & 
	+  \frac{\tems}{\ep^2}\sigma_\ep \left[[0,t]\times \overline\Omega_\ep\right]\\
	& = \int_{\Omega_\ep} \left( \frac{1}{2}\vrez|\uez|^2 \right) \dx 
	\\ &+ \frac{1}{\ep^2}\left(H_{\tems}(\vrez,\temz) - (\vrez - \vret)\frac{\partial H_{\tems}(\vret,\tems)}{\partial \vr} 
	- H_{\tems}(\vret,\tems)\right) \dx,
	\end{split}
	\end{equation} 
where $H_{\tems}$ is ballistic free energy and 
	$
	H_{\tems}(\vr,\temp)=\vr \left( e(\vr,\temp) - \tems s(\vr,\temp)  \right).
	$

It is provided (see Lemma 5.1 in \cite{FN}) that 
	$H_{\tems}(\vre,\tem) - (\vre - \vret)\frac{\partial H_{\tems}(\vret,\tems)}{\partial \vr} 
	- H_{\tems}(\vret,\tems)$
is non-negative, strictly coercive, attain global minimum zero at point $(\vret,\tems)$, 
dominates internal energy $\vr e$ and entropy $s$ far from $(\vret,\tems)$.
Therefore according to our assumptions
we are able to deduce  from \eqref{H2} that (for details see \cite{FN,FS,FKKS,AWK})
	\begin{equation*}\label{est1}
	\ess \sup\limits_{t\in (0,T)} \int_{\Omega_\ep} \vre |\ue|^2(t,\cdot) \dx \leq c, \quad \ess \sup\limits_{t\in(0,T)} \| \sqrt{\vre}\ue\|_{L^2(\Omega_\ep;\R^3)} \leq c 
	\end{equation*}
	\begin{equation*}\label{est2}
	\ess \sup\limits_{t\in (0,T)} \left\| \left[ \frac{\vre - \vret}{\ep}\right]_\ess (t,\cdot) \right\|_{L^2(\Omega_\ep)} \leq c ,
\quad
	\ess \sup\limits_{t\in (0,T)} \left\| \left[ \frac{\tem - \tems}{\ep}\right]_\ess (t,\cdot) \right\|_{L^2(\Omega_\ep)} \leq c ,
	\end{equation*}
	\begin{equation*}
	\| \sigma_\ep\|_{{\mathcal{M}}^+ ([0,T]\times\Omega_\ep )} \leq \ep^2 c ,
	\end{equation*}
	\begin{equation*}\label{est5}
	\ess \sup\limits_{t\in (0,T)} \int_{\Omega_\ep}
	\left(
	| [ \vre e(\vre,\tem)]_\res| +  | [ p(\vre,\tem)]_\res|   + | [ \vre s(\vre,\tem)]_\res| \dx
	\right) 
	 \leq  \ep^2 c ,
	\end{equation*}
	\begin{equation*}\label{est31}
	\ess \sup\limits_{t\in (0,T)} \int_{\Omega_\ep} [ \vre]^{5/3}_\res (t,\cdot)  +  [ \tem]^{4}_\res (t,\cdot) \dx \leq \ep^2 c,
	\quad
	\ess \sup\limits_{t\in (0,T)} \int_{\Omega_\ep}	\bbbone_\res(t, \cdot)  dx \leq \ep^2 c,  
	\end{equation*}
	\begin{equation*}
		\ess \sup\limits_{t\in (0,T)} \left\| \left[  \frac{\vre - \vret}{\ep}   \right]_\res  \right\|_{L^1(\Omega_\ep)} \leq c \ep ,
 	\end{equation*} 
	\begin{equation*}\label{est55}
	\int_0^T \left\| \frac{\tem - \tems}{\ep} \right\|^2_{W^{1,2}(\Omega_\ep;R^3)} \dt +
	\int_0^T \left\| \frac{\log(\tem) - \log(\tems)}{\ep} \right\|^2_{W^{1,2}(\Omega_\ep;R^3)} \dt
	< c,
	\end{equation*}
	\begin{equation*}\label{uw12}
	\int_0^T \left\| \ue \right\|^2_{W^{1,2}(\Omega_\ep; \R^3)} \dt < c. 
	\end{equation*}
 
 \subsection{Convergence}
 The hypotheses stated on the family of $\{ \Omega_\ep \}_\ep$ provides us: 
 \begin{itemize} 
 \item the uniform extension property \cite{Jones}. Namely there exists an extension operator $E_\ep$ s.t. $E_\ep: W^{1,p}(\Omega_\ep) \mapsto W^{1,p}(\R^3)$, 
 $E_\ep [v]|_{\Omega_\ep} = v$ and $\| E_\ep [v] \|_{W^{1,p}(\R^3)} \leq  c \| v \|_{W^{1,p}(\Omega_\ep)}$, where the constant $c$ is independent of $\ep \to 0$.
 \item there exists bounded domain $O$ s.t. $\R^3 \setminus O$ satisfy  the uniform $\alpha$-cone condition 
and  a suitable subsequence of $\ep$'s  such that  $| (\R^3 \setminus O_\ep) \setminus (\R^3 \setminus O)| \to 0$  as $\ep \to 0 $. This property is crucial when studying stability if the spectral properties of the Neumann Laplacian, see \cite{FKKS,AK}, to provide decay of acoustic waves. For each $x_0 \in \partial  O$ there is $x_{\ep,0} \in \partial O_\ep$ such that $x_{\ep,0} \to x_{0}$ and $O \subset B_s(0)$ and for any compact $K \subset \Omega$, there exists $\ep(K)$ such that $K \subset \Omega_\ep$ for all $\ep < \ep(K).$
\end{itemize}
Since the family of $\{ \Omega_\ep \}_\ep$ possesses a uniform extension property we may deduce from uniform estimates that
	\begin{equation}\label{weakU}
	\ue \weak \vec{U} \quad \mbox{ weakly in } L^2(0,T;W^{1,2}(\R^3;\R^3)),
	\end{equation}
	\begin{equation*}
	\ess \sup\limits_{t\in(0,T)} \| \tem(t,\cdot) - \tems \|_{L^2(\Omega_\ep) } \to 0 \mbox{ as } \ep \to 0,
	\end{equation*}
	\begin{equation*}
	\Theta_\ep = \frac{\tem - \tems}{\ep} \weak \Theta \mbox{ weakly in } L^2(0,T;W^{1,2}(\R^3)) 	
	\end{equation*}
Following the same procedure as in \cite{FS,AWK}  by uniform estimates and closeness of $\vrs$ and $\vret$ \eqref{nff} 
we get 	
		\begin{equation*}\label{ro531}
		\begin{split}
	& \ess \sup\limits_{t\in(0,T)} \| \vret(t,\cdot)-\vrs \|_{L^{5/3}+ L^q (\Omega_\ep)}\to 0, 
	\\ & 
	\ess \sup\limits_{t\in(0,T)} \| \vre(t,\cdot) - \bar\vr \|_{L^2 + L^{5/3 }+ L^q (\Omega_\ep)} \to 0 \mbox{ as } \ep\to 0
	 \quad \mbox{ for } q>3,
	 \end{split}
	\end{equation*}	
	\begin{equation*}\label{k1}
	\frac{\vre - \vrs}{\ep} \weakstar r \mbox{ weakly* in } L^\infty(0,T; L^{5/3}(K)) \mbox{ for  any compact } K\subset \Omega \mbox{ for } \ep \to 0.
	\end{equation*} 
Therefore fluid density becomes constant since $\ep\to 0$, i.e. as the Mach number tends to zero.
Then continuity equations provides, that  
	$$\div \U = 0 \mbox{ a.a. in } (0,T) \times \Omega .$$
By boundary conditions and properties of $\Omega_\ep$ the limit velocity field satisfies the impermeability condition $\vec{U} \cdot \vec{n}|_{\partial O} = 0$ in a weak sense. 
Moreover the analysis provided by \cite{BFNW}, see also \cite[Sec. 6.2]{FKKS}, gives that $\vU|_{\partial O} = 0$ if (D6) is satisfied. 

To pass to the limit in rescaled NSF$_\ep$ system one of the most difficult steps is to provide strong convergence of   
the velocity field in order to control the limit of convective term. Namely we need to show that
		\begin{equation*}\label{stue21}
		\ue \to \vU \quad \mbox{strongly in } L^2((0,T)\times K) \mbox{ for any compact } K\subset \R^3 \setminus O. 
		\end{equation*}
The main obstacle here are possible oscillations in time of the momentum, since from momentum equations we do not control its time derivative. Then one can observe that it is sufficient to provide that (see \cite{F_LDB,FS,AWK})
	\begin{equation}\label{stue22}
	\vre \ue \to \vrs\vU \quad\mbox{in } L^2(0,T; W^{-1,2}(K)).
	\end{equation}
Then due to \eqref{weakU} it is even enough to prove, instead of  \eqref{stue22}, that
	$$\{  t \to \int_{\R^3} (\vre \ue) (t, \cdot) \vec\varphi \dx\} \mbox{ is precompact in }L^2(0,T)$$ 
and 
	\begin{equation}\label{c1}
	\left\{ t \to \int_{\R^3} \vre \ue (\cdot, t) \cdot \vec\varphi \dx \right\} \quad 
	\to \quad 
	\left\{ t \to \vrs \int_{\R^3} \vec{U}(\cdot, t) \cdot \vec\varphi \dx \right\}
	\mbox{ in }L^2(0,T)
	\end{equation}
for any fixed $\vec\varphi \in C^\infty(\R^3)$ where ${\rm supp\,} \vec\varphi \subset K$ as $\ep \to 0.$ 

\subsection{Reformulation to the wave equation. Dispersive estimates - local decay of acoustic wave}

As it was already emphasised, our aim now is to show \eqref{c1}. 
This will be provided  by the analyse of Lighthill's acoustic analog (see \cite{Lighthill}) of our primitive NSF$_\ep$ system, namely 
	\begin{equation}\label{ac2}
	\ep\partial_t  S_\ep + \omega \div \vec{V}_\ep  = \ep \tilde{f}_\ep^1, \quad \quad\ep\partial_t \vec{V}_\ep + \nabla_x S_\ep         = \ep \vec{\tilde{f}}_\ep^2,
	\end{equation}
with homogenous  Neuman boundary condition
	$ \vec{V}_\ep \cdot \n |_{\partial \Omega_\ep} = 0$
 where
	\begin{equation*}\label{sep}
	S_\ep 
	= A \left( \frac{\vre - \vrs}{\ep}\right) + B \left(  \frac{ \vre s(\vre,\tem) - \vrs s(\vrs,\tems) }{ \ep } \right)
	-\vrs F_\ep  + \frac{B}{\ep} \Sigma_\ep ,
	\quad \quad
	\vec{V}_\ep = \vre\ue ,
	\end{equation*}
	\begin{equation*}\label{ff1e2}
	\tilde{f}_\ep^1 =  
	\div 
	 \underbrace{  B
	\left( 
	\vre \frac{s(\vrs,\tems) - s(\vre,\tem)}{\ep}\ue 
	\right) 
	}_{H_\ep^1}
	+ \div 
	 \underbrace{ B
	\left( \frac{\kappa(\tem)}{\tem} \frac{\nabla_x \tem}{\ep} \right)
	}_{H_\ep^2} 
	\end{equation*}
	\begin{equation*}\label{ff2e2}
	\begin{split}
	{\vec{\tilde f}}^2_\ep = &
	\nabla_x
	\underbrace{
	  \frac{1}{\ep} 
	 \left[
	 A \left( \frac{\vre - \vrs}{\ep}\right) + B \left(  \frac{ \vre s(\vre,\tem) - \vrs s(\vrs,\tems) }{ \ep } \right) 
	 - \left(  \frac{p(\vre,\tem) - p(\vrs,\tems)}{\ep} \right)
	\right] 
	}_{G_\ep^3}
	 \\ & - \div \left[
	 \underbrace{
	 (\vre\ue \otimes \ue )
	 }_{G_\ep^{2,2}} +
	 \underbrace{\tS_\ep }_{G_\ep^{2,1}} \right]
	 + 
	 \underbrace{ \frac{\vre - \vrs}{\ep}\nabla_x F_\ep }_{G_\ep^4}
	 + 
	 \underbrace{ B \frac{1}{\ep^2}  \nabla_x \Sigma_\ep }_{\nabla_x G_\ep^1}.
	\end{split}
	\end{equation*}
Where $\Sigma_\ep$ is a time lifting of $\sigma_\ep$ (\cite{FN,FS,AWK}) and constants $A$, $B$, $\omega$ are chosen s.t. 
	$
	B\vrs \partial_\temp s(\vrs,\tems) = \partial_\temp p(\vrs,\tems)$
	{and} $
	A + B \partial_\vr (\vr s )(\vrs,\tems) = \partial_\vr p  (\vrs,\tems),$ $ \omega = \partial_\vr p (\vrs, \tems) + 
	\frac{ | \partial_\temp p (\vrs,\tems) |^2}{ \vrs^2 \partial_\temp s (\vrs,\tems) }  
	> 0$
(see e.g.  \cite{FN,FS,AWK}). Notice that $\omega$ is bounded due to structural restrictions on $p$ and $s$.

Let 
$\nabla_x \Phi_\ep$ denote acoustic potential, i.e. 
	$$
	\vec{V}_\ep = \vec{H}_\ep[\vec{V}_\ep] + \nabla_x \Phi_\ep.
	$$
Accordingly  we may  rewrite  (\ref{ac2})$_1$ in the following form
 \begin{equation}\label{aceqw12}
 	\begin{split}
	\ep \int_0^T \langle S_\ep(t,\cdot),\partial_t \varphi \rangle \dt  
	 & + \omega \int_0^T \int_{\Omega_\ep} \nabla_x \Phi_\ep \cdot \nabla_x\varphi \dxdt 
	\\ & =
	 \ep \ \langle S_{0,\ep}, \varphi(0,\cdot) \rangle 
	+ 
	\ep  \int_0^T \int_{\Omega_\ep}
	( H^1_\ep  +  H^2_\ep ) \cdot \nabla_x \varphi \dxdt 
	\end{split}
	\end{equation}
	for all $\varphi \in C^\infty_c ( [0,T] \times \overline{\Omega}_\ep )$. 
Next since $\varphi = \nabla_x \len^{-1}[\varphi]$ is an admissible test function in \eqref{ac2}$_2$
(due to slip boundary condtion on $\ue$)
 we obtain by integration by parts that 
	\begin{equation}\label{aceqw22}
	\begin{split}
	&
	\ep \int_0^T \int_{\Omega_\ep} \Phi_\ep \cdot \partial_t \varphi  \dt  
	-  \int_0^T  \langle S_\ep ,  \varphi  \rangle _{[{{\mathcal M}}, C]} \dt 
	= -  \ep \int_{\Omega_\ep} V_{0,\ep} \cdot \nabla_x \len^{-1}[\varphi(0,\cdot)] \dx 
	\\ & 
	- \ep \Big\{
	 \int_0^T \langle G^1_\ep (t, \cdot),  \varphi   \rangle   \dt +
	 \int_0^T \int_{\Omega_\ep}
	G^{2,1}_\ep : \nabla_x^2 \len^{-1} [ \varphi] 
	\dxdt 
	\\ &
	+ \int_0^T \int_{\Omega_\ep}
	G^{2,2}_\ep : \nabla_x^2 \len^{-1} [ \varphi] \dxdt 
	+  \int_0^T \int_{\Omega_\ep} G^3_\ep  \varphi \dxdt
	\\& 
	+  \int_0^T \int_{\Omega_\ep} G^4_\ep \cdot \nabla_x \len^{-1} [ \varphi] \dxdt \Big\} .
	\end{split}
	\end{equation}
The above equations represent a weak formulation of the acoustic equation for the potential of the gradient part of the momentum with Neumann boundary condition. 
\smallskip

Summarising computation from previous sections, due to uniform estimates obtained in Section~\ref{unib} 
equations \eqref{aceqw12} and \eqref{aceqw22} can be rewritten in the following more conscious form (see \cite{FKKS,AWK})
	 \begin{equation}\label{aceqw13}
	\begin{split}
	&
	\ep \int_0^T \langle S_\ep(t,\cdot),\partial_t \varphi \rangle \dt  
	+ \omega \int_0^T \int_{\Omega_\ep} \nabla_x \Phi_\ep \cdot \nabla_x\varphi \dxdt \\
	& 
	= \ep \ \langle S_{0,\ep}, \varphi(0,\cdot) \rangle 
	+ \frac{\ep}{\ep^{2\beta}} 
	  \int_0^T \int_{\Omega_\ep} J^{1}_\ep \varphi 
	+ J^{2}_\ep (-\len)^{3/2}[\varphi]
	\\ & + J^{3}_\ep (-\len)^{1/2}[\varphi]  +
	 J^{4}_\ep (-\len)[\varphi] \dxdt
	\end{split}
	\end{equation}
	for all $\varphi \in C^\infty_c ( [0,T] \times \overline{\Omega}_\ep )$ and 
	\begin{equation}\label{aceqw23}
	\begin{split}
	&
	\ep \int_0^T \int_{\Omega_\ep} \Phi_\ep \cdot \partial_t \varphi  \dt  
	-  \int_0^T  \langle S_\ep ,  \varphi  \rangle \dt 
	= - \ep \int_{\Omega_\ep} \Phi_{0,\ep} \varphi(0,\cdot) \dx 
	\\
	&
	- \frac{\ep}{ \ep^{2\beta} } \int_0^T \int_{\Omega_\ep} 
	\Big\{
	\tilde{J}^1_\ep \varphi 
	+ \tilde{J}^2_\ep (-\len)^{-1/2}[\varphi]
	+ \tilde{J}^3_\ep (-\len)^{1/2}[\varphi] 
	\\ & +\tilde{J}^4_\ep (-\len)^{-1}[\varphi]
	 + \tilde{J}^5_\ep (-\len)[\varphi]
	 \Big\}
	  \dxdt
	 \end{split}
	 \end{equation}
for any $\varphi \in C^\infty_c ([0,T) \times K)$, $K$ compact subset of $\R^3 \setminus O$, 
$\nabla_x \varphi \cdot \n |_{\partial \Omega_\ep} =0$,
where 
	$$
	\| {J}^{i} \|_{L^2( (0,T) \times \Omega_\ep )} < c \mbox{ for } i=1,\dots, 4
	\mbox{ and } 
	\| \tilde{J}^j \|_{L^2( (0,T) \times \Omega_\ep )} < c
	\mbox{ for } j=1,\dots, 5
	$$
and for sufficiently small $\ep$ and supplemented with the following initial data	$$S_{0,\ep} = (-\len)[\tilde{S}_{0,\ep}^1] + (-\len)^{1/2} [\tilde{S}_{0,\ep}^2] + \tilde{S}_{0,\ep}^3,$$  with  
$\|\tilde{S}_{0,\ep}^i \|_{L^2(\Omega_\ep)} \leq c$  and $$\Phi_{0,\ep} = (-\len)^{-1} \div V_{0,\ep}, \mbox{ where }\| (-\len)^{-1/2} [\Phi_{0,\ep}] \|_{L^2(\Omega_\ep)} \leq c.$$

Then the Duhamel formula gives as an explicit formulation for acoustic potential, i.e.: 
	\begin{equation}\label{duhamel1}
	\begin{split}
	\Phi_\ep (t,\cdot) & = \frac{1}{2} \exp \left(  \pm i\sqrt{- \omega \len} \frac{t}{\ep} \right)
	\left[ \Phi_{0,\ep} \pm  \frac{ i }{ \sqrt{- \omega \len} }[ S_{0,\ep} ]\right] 
	\\ &
	+ \ep^{-2\beta}  \frac{1}{2} \int_0^T \exp \left( \pm  i\sqrt{- \omega \len} \frac{t-s}{\ep} \right)
	\left[ \tilde{F}_{2,\ep}(s) \pm \frac{ i }{ \sqrt{- \omega \len} } \tilde{F}_{1,\ep}(s) \right] {\rm d} s,
	\end{split}
	\end{equation}
where 	
	\begin{equation*}
	\begin{split}
	\tilde{F}_{1,\ep} &= J^{1}_\ep  
	+ (-\len)^{3/2} [J^{2}_\ep ]
	+ (-\len)^{1/2} [J^{3}_\ep]  +
	(-\len) [J^{4}_\ep] , 
	\\
	\tilde{F}_{2,\ep}  & =  \tilde{J}^1_\ep  
	+ (-\len)^{-1/2} [\tilde{J}^2_\ep] 
	+ (-\len)^{1/2} [\tilde{J}^3_\ep] +
	 (-\len)^{-1} [\tilde{J}^4_\ep] +   (-\len) [\tilde{J}^5_\ep]
	 \end{split}
	 \end{equation*}
(see \eqref{aceqw13}, \eqref{aceqw23}). Let us remark that the ''large'' coefficient $\ep^{- 2\beta}$ appearing in \eqref{aceqw13}, \eqref{aceqw23} and \eqref{duhamel1} is a consequence or roughness of the obstacle $O_\ep$ (see (D5)). More precisely, an elliptic estimate employed to derive  \eqref{aceqw13}, \eqref{aceqw23} depends on $\ep$, i.e.
	$
	\|  \nabla_x^2 \varphi \|_{L^p(\Omega_\ep)}\leq c(p) \left(  \| \Delta_x \varphi \|_{L^p({\Omega}_\ep)}  + 
	 \frac{1}{\ep^{2\beta}}\| \varphi \|_{L^p({\Omega}_\ep)}  \right)
	$
for any $\varphi \in C^\infty_c(\overline{\Omega_\ep})$ with $\nabla_x \varphi \cdot \n|_{\partial {\Omega}_\ep} = 0$, with $1<p<\infty .$

With above formulation at hand and by methods developed in \cite{F_Ld} we are able to provide local decay of acoustic wave and consequently  to show that 
	\begin{equation}\label{pg0}
	\left\{ t \to \int_{\Omega_\ep} \Phi_\ep G(- \len ) [ \varphi ] \dx   \right\} \to 0 \quad \mbox{in } L^2(0,T),
	\end{equation}
any $G\in C^\infty_c(0,\infty)$, what in fact  is a key point to prove \eqref{c1} and consequently to provide convergence in convective term (see for details \cite{FKKS,AWK}).
The following lemma gives a local decay of acoustic waves.
\begin{lemma}[\cite{F_Ld,FKKS}]\label{de}
We have 
	\begin{equation*}\label{de1}
	\int_0^T \left| \left\langle  \exp\left(  i \sqrt{-\len} \frac{t}{\ep} [\Psi] , G(- \len ) [\varphi] 
	\right) \right\rangle_{\widetilde\Omega_\ep}  \right|^2 \dt \leq \ep c(\varphi, G) \| \Psi \|^2_{L^2(\widetilde\Omega_\ep)}
	\end{equation*}
for any $\varphi \in C^\infty_c(K)$, $\Psi \in L^2(\tilde\Omega_\ep)$, and any $G\in C^\infty_c(0,\infty)$, where is s.t. $\overline{K} \subset \R^3 \setminus O_\ep$.
\end{lemma}

Lemma~\ref{de} applied to $\Phi_\ep$ given by formula \eqref{duhamel1} provides \eqref{pg0}, if $\beta < \frac{1}{4}$, see \cite{F_Ld} for details.
The explicitly given rate of the decay in Lemma~\ref{de} allow to compensate 
exploding coefficient $\ep^{-2\beta}$ which reflects the influence of perturbations of the domain.
Moreover, let us remark that in order to provide good properties of the spectrum of Neumann Laplacian $-\len$ it is crucial to notice that the outer boundary (the boundary of the sphere $\mathcal{S}_\ep$) is irrelevant for the local analysis (on supports of test functions $\vec{\varphi}$)  and in fact we may consider the operator  $-\len$ on unbounded domain $\R^3 \setminus O_\ep$. Indeed in \eqref{ac2} the speed of propagation is finite and proportional to $\sqrt{\omega} / \ep $ and the boundary $\mathcal{S}_\ep$ is sufficiently ''far'', since $\delta>1.$
For details see again \cite{FKKS,AWK}.

\def\ocirc#1{\ifmmode\setbox0=\hbox{$#1$}\dimen0=\ht0 \advance\dimen0
  by1pt\rlap{\hbox to\wd0{\hss\raise\dimen0
  \hbox{\hskip.2em$\scriptscriptstyle\circ$}\hss}}#1\else {\accent"17 #1}\fi}

\end{document}